\documentclass [12pt,fleqn]{article}
\usepackage{latexsym,amsmath}
\title{Numerical solutions of
integrodifferential systems by hybrid of general block-pulse
functions and the second Chebyshev polynomials}

\author{Xing Tao Wang*\\
\noindent\footnotesize{\it Department of Mathematics, Harbin
Institute of Technology, Harbin $150001$, P.R. China}}

\date{}

\begin{document}
\footnotetext{\hspace*{-17pt}*Corresponding author. \ Tel:
+86-0451-86413074; Fax: +86-0451-86414216
\\ \textit{E-mail address}: xingtao@hit.edu.cn (Xing Tao Wang)
} \maketitle \noindent\rule[4pt]{6.02in}{0.5pt}
{{\small\textbf{Abstract } By applying hybrid functions of general
block-pulse functions and the second Chebyshev polynomials,
 integrodifferential systems  are
converted into a system of algebraic equations. The approximate
solutions of  integrodifferential systems are derived. The
numerical examples illustrate that the algorithms are valid.

\noindent\rule[4pt]{6.02in}{0.5pt}
\\
\\

\noindent\textbf{1. \ Introduction}
\\

Integrodifferential  systems are a class of systems of importance
\cite{Kochetkov}. Since it is  difficult to obtain the analytic
solutions of integrodifferential  systems, numerical methods to
obtain  approximate solutions are interested. Many orthogonal
functions or polynomials, such as block-pulse functions
\cite{block1,block2}, Walsh functions \cite{Walsh}, Fourier series
\cite{Fourier} , Legendre polynomials \cite{Legendre}, Chebyshev
polynomials \cite{Chebyshev} and Laguerre polynomials
\cite{Laguerre}, were used to derive solutions of some systems. In
recent years  different kinds of hybrid functions
\cite{hybridTaylor,hybridLegendre} and Legendre wavelets
\cite{Legendrewavelets} were used. In this article we apply the
hybrid functions consisting of general block-pulse functions and
the second Chebyshev polynomials to converting
 integrodifferential systems into a system of algebraic equations. The hybrid solutions of
integrodifferential systems are obtained. We present the general
operational matrices. The numerical examples show that this method
is convenient for application.
\\
\\

\noindent\textbf{2. Preliminaries}
\\

\noindent\textit{2.1. Definitions}
\\

A set of block-pulse function $b_k (t), \ k=1,2,\dots ,K$, on the
interval $[t_0, t_f )$ are defined as
\[
b_k (t) = \bigg\{ \begin{array}{l}
 1, \ t_{k-1} \le t < t_k , \  \\
 0, \ {\rm otherwise},
 \end{array}
\]
where  $t_K=t_f$ and $[t_{k-1},  t_k)\subset[t_0, t_f )$, $k =
1,2, \dots ,K$.

The $m$ order second Chebyshev polynomials  in the interval $[ -
1,1]$ are defined by the following:
\begin{equation}
S_m (t) = {\rm sin}((m+1) {\rm arccos}(t))(1-t^2)^{-\frac{1}{2}},
\end{equation}
which are also given by the following recursive formula
\begin{equation}
\bigg\{\begin{array}{l}
S_0 (t) = 1, \ S_1 (t) = 2t, \\
S_{m + 1} (t) = 2tS_m (t) - S_{m - 1} (t), \ m = 1,2, \dots
.\label{chebyshev}
\end{array}
\end{equation}

 The hybrid functions $h_{km}
(t), \ k = 1,2, \dots ,K; \ m = 0,1, \dots ,M - 1$, on the
interval $[t_0, t_f )$ are defined as
\begin{equation}
h_{km} (t) = b_k (t)S_m (d_k^{-1}( 2t - t_{k-1}- t_k)),
\end{equation}
where $d_k=t_k-t_{k-1}$, $k = 1,2, \dots ,K$.
\\

\noindent\textit{2.2. The operational matrix}
\\

Let
\[
H_k (t) = [h_{k0} (t), \dots ,h_{k,M - 1} (t)]^\tau ,\quad H(t) =
[H_1^\tau (t), \dots ,H_K^\tau (t)]^\tau.
\]
It is not difficult to  verify that
\begin{equation}
\int_{\,t_0}^{\,t} {H(\bar{t}){\rm  d}\bar{t} \simeq PH(t)} ,
\end{equation}
where
\[
 P=\frac{1}{2}{\rm diag}(d_1,\dots,d_K) \otimes
\hat{P}+\sum\limits_{k = 1}^{K - 1} {\sum\limits_{i = 1}^{K - k}
d_i{E_{i,i + k}^{(K)} } } \otimes \sum\limits_{k =
1}^{[\frac{M+1}{2} ]} \frac{1}{2k-1}E_{2k-1,1}^{(M)},
\]
\[
\hat{P}=E_{11}^{(M)} -\frac{3}{4}E_{21}^{(M)}
+\frac{1}{2}\bigg(\sum\limits_{k = 1}^{M - 1} \frac{1}{k}E_{k,k +
1}^{(M)} -\sum\limits_{k = 2}^{M - 1}  \frac{1}{k +1}E_{k +
1,k}^{(M)} \bigg)+\sum\limits_{k = 3}^{M }
\frac{(-1)^{k-1}}{k}E_{k1}^{(M)},
\]
$E_{ij}^{(m)} $ is the $m\times m$ matrix with 1 at its entry
$(i,j)$ and zeros elsewhere, $\otimes$ denotes Kronecker product
and $[\frac{M+1}{2} ]$ is the greatest integer part of
$\frac{M+1}{2} $.
\\

\noindent\textit{2.3. Function approximation}
\\

An $l$-dimensional vector  function $f(t)$ on the interval $[t_0,
t_f )$ is expressed as
\begin{equation}
f(t) \simeq \sum\limits_{k = 1}^K {\sum\limits_{m = 0}^{M - 1}
{f_{km} h_{km} (t)} },\label{f}
\end{equation}
where
\begin{equation}
f_{km} = \frac{2}{\pi }\int_{\,{-1}}^{\,1 } {f(2^{-1}(d_k t +
t_{k-1}+ t_k))S_m (t)(1-t^2)^{\frac{1}{2}}\text dt}.
\end{equation}
Rewrite $f(t)$ as
\[
f(t) \simeq \sum\limits_{k = 1}^K {F_k H_k (t) = } FH(t),
\]
where
\[
F_k = [f_{k0} , \dots ,f_{k,M - 1} ], \quad  F = [F_1 , \dots ,F_K
].
\]
For corresponding $F_k $ and $F$ we denote
\[
\hat {F}_k = [f_{k0}^\tau , \dots ,f_{k,M - 1}^\tau ]^\tau, \quad
\hat {F} = [\hat {F}_1^\tau , \dots ,\hat {F}_K^\tau ]^\tau,
\]
where $\tau$ is the transpose.
\\

\noindent\textit{2.4. Expression of the product of a matrix
function and a vector function}
\\

Let a matrix function $M(t)$ be appropriate to a vector function
$f(t)$. We express $M(t)$ and $f(t)$, respectively, as
\[
M(t) \simeq \sum\limits_{k = 1}^K {\sum\limits_{m = 0}^{M - 1}
{M_{km} h_{km} (t)} } , \quad f(t) \simeq \sum\limits_{k = 1}^K
{\sum\limits_{m = 0}^{M - 1} {f_{km} h_{km} (t)} }.
\]
Then
\[
M(t)f(t) \simeq \sum\limits_{k = 1}^K {\sum\limits_{i = 0}^{M - 1}
{\sum\limits_{j = 0}^{M - 1} {M_{ki} f_{kj} h_{ki} (t)} } } h_{kj}
(t).
\]
From
\[
h_{ki} (t)h_{kj} (t) \simeq \sum\limits_{m = 0}^{M - 1}
{d_{km}^{(ij)} h_{km} (t)},
\]
where
\begin{equation}
d_{km}^{(ij)} = \frac{2}{\pi }\int_{\,-1}^{\,1} {S_i (t)S_j (t)S_m
(t)(1-t^2)^{\frac{1}{2}}\text dt},
\end{equation}
we have
\begin{gather}
\begin{split}
M(t)f(t) \simeq \sum\limits_{k = 1}^K {\sum\limits_{m = 0}^{M - 1}
{\bigg(\sum\limits_{i = 0}^{M - 1} {\sum\limits_{j = 0}^{M - 1}
{d_{km}^{(ij)} M_{ki} f_{kj}\bigg)h_{km} (t)} } } }\hspace*{-15pt}\\
= \sum\limits_{k = 1}^K {\sum\limits_{m = 0}^{M - 1} {\tilde
{M}_{km} h_{km} (t)} } = \sum\limits_{k = 1}^K {\tilde {M}_k } H_k
(t),\label{Mf}
\end{split}
\end{gather}
where
\[
\tilde {M}_{km} = \sum\limits_{i = 0}^{M - 1} {\sum\limits_{j =
0}^{M - 1} {d_{km}^{(ij)} M_{ki} f_{kj} } } = \hat {M}_{km} \hat
{F}_k,\quad\hat {F}_k = [f_{k0}^\tau, \dots, f_{k,M - 1}^\tau
]^\tau, \ \hat {\tilde {M}}_k = \hat {M}_k \hat{F}_k,
\]
\[
\hat {\tilde {M}}_k = [\tilde {M}_{k0}^\tau , \dots ,\tilde
{M}_{k,M - 1}^\tau ]^\tau , \quad \hat {M}_k = [\hat {M}_{k0}^\tau
, \dots ,\hat {M}_{k,M - 1}^\tau ]^\tau,
\]
\[
\hat {M}_{km} =\bigg[\sum\limits_{i = 0}^{M - 1} {d_{km}^{(i0)}
M_{ki} } , \dots ,\sum\limits_{i = 0}^{M - 1} {d_{km}^{(i,M - 1)}
M_{ki} }\bigg].
\]
Therefore
\begin{gather}
N(\bar{t},t)f(t) \simeq \sum\limits_{k = 1}^K {\sum\limits_{m =
0}^{M - 1} {\tilde {N}_{km}(\bar{t}) h_{km} (t)} }\label{Nf},
\end{gather}
where
\[
\tilde {N}_{km} (\bar{t}) = \hat {N}_{km}(\bar{t}) \hat
{F}_k,\quad\hat {F}_k = [f_{k0}^\tau, \dots, f_{k,M - 1}^\tau
]^\tau,
\]
\[
\hat {N}_{km}(\bar{t}) =\bigg[\sum\limits_{i = 0}^{M - 1}
{d_{km}^{(i0)} N_{ki}(\bar{t}) } , \dots ,\sum\limits_{i = 0}^{M -
1} {d_{km}^{(i,M - 1)} N_{ki}(\bar{t}) }\bigg],
\]
\[
N_{ki}(\bar{t})  = \frac{2}{\pi }\int_{\,{-1}}^{\,1 }
{N(\bar{t},2^{-1}(d_k t + t_{k-1}+ t_k))S_i
(t)(1-t^2)^{\frac{1}{2}}\text dt}.
\]
Let
\begin{equation}
w(t)  = \int_{\,{t_0}}^{\,t_f } {N(t,\bar{t})f (\bar{t})\text
d\bar{t}} \simeq \sum\limits_{k = 1}^K {\sum\limits_{m = 0}^{M -
1} {w_{km} h_{km} (t)} },\label{Nx}
\end{equation}
\[
\hat {N}_{km}(t)  \simeq \sum\limits_{k = 1}^K {\sum\limits_{m =
0}^{M - 1} {\hat {N}^{(jl)}_{km} h_{jl} (t)} },
\]
where
\[
\hat {N}^{(jl)}_{km} =\bigg[\sum\limits_{i = 0}^{M - 1}
{d_{km}^{(i0)} N^{(jl)}_{ki} } , \dots ,\sum\limits_{i = 0}^{M -
1} {d_{km}^{(i,M - 1)} N^{(jl)}_{ki} }\bigg],
\]
\[
N^{(jl)}_{ki}  = \frac{2}{\pi }\int_{\,{-1}}^{\,1 }
{N_{ki}(2^{-1}(d_j t + t_{j-1}+ t_j))S_l
(t)(1-t^2)^{\frac{1}{2}}\text dt}.
\]
By Eq.(\ref{Nf}) we have
\begin{equation}
w(t)    \simeq \sum\limits_{j = 1}^K {\sum\limits_{l = 0}^{M -
1}\sum\limits_{k = 1}^K \sum\limits_{m = 0}^{M - 1} {\hat
{N}^{(jl)}_{km}\bigg({\int_{\,{t_{k-1}}}^{\,t_k }h_{jl} (t){\rm
d}t}\bigg)\hat{F}_k h_{jl} (t)} }.
\end{equation}
So
\begin{equation}
\hat{W}  = \sum\limits_{k = 1}^{[\frac{M+1}{2} ]} \sum\limits_{i =
1}^K \sum\limits_{j = 1}^{K}E^{(K)}_{ij}\otimes
\frac{d_j}{2k-1}{{N}^{(i)}_{j,2k-1}\hat{F}},\label{W}
\end{equation}
where
\[
{N}^{(i)}_{j,2k-1}=\bigg[{{N}^{(i0)}_{j,2k-1}}\hspace{-10pt}^\tau\hspace{10pt},\dots,{{N}^{(i,M-1)}_{j,2k-1}}^\tau\bigg]^\tau.
\]
\\

\noindent\textbf{3. Analysis of integrodifferential systems }
\\

Consider the following integrodifferential system \cite{Kochetkov}
described by
\begin{equation}
\Bigg\{ \begin{array}{l}
 \dot {x}(t) = A(t)x(t) + \displaystyle\int_{\,{t_{0}}}^{\,t_f }{N (t,\bar{t})x(\bar{t}){\rm d}\bar{t}}+B(t)u(t), \
  t \in [t_0, t_f ],\\
 x(t_0) =x_0,\label{sys}
 \end{array}
\end{equation}
where $x(t)$ is an {\it n}-dimensional vector function and $u(t)$
an $r$-dimensional vector function. $ A(t)$, $B(t)$ and $
N(t,\bar{t})$ are the matrix functions of appropriate dimensions.
Express $ A(t)$, $x(t)$, $\int_{\,{t_{0}}}^{\,t_f }N
(t,\bar{t})x(\bar{t}){\rm d}\bar{t}$, $ B(t)$ and $u(t)$
respectively, as Eqs.(\ref{f}) and (\ref{Nx}). By Eqs.(\ref{Mf})
and (\ref{Nf}) we have
\[
A(t)x(t) \simeq \sum\limits_{k= 1}^K {\sum\limits_{m = 0}^{M - 1}
{\tilde {A}_{km} h_{km} (t)} }, \quad B(t)u(t) \simeq
\sum\limits_{k= 1}^K {\sum\limits_{m = 0}^{M - 1} {\tilde {B}_{km}
h_{km} (t)} } ,
\]
\[
\int_{\,{t_0}}^{\,t_f } {N(t,\bar{t})f (\bar{t}){\rm d}\bar{t}}
\simeq \sum\limits_{k = 1}^K {\sum\limits_{m = 0}^{M - 1} {w_{km}
h_{km} (t)} }.
\]
Integrating Eq.(\ref{sys}) from $t_0$ to $t$ and combining
Eqs.(\ref{Mf}), (\ref{Nf}) and (\ref{W})
 we obtain
\[
[X_1 , \dots ,X_K ]H(t) - [X_{01} , \dots ,X_{0K} ]H(t)
\]
\[
 \hspace*{20pt}= \int_{\,t_0}^{\,t} {\{[  \tilde {A}_1 , \dots ,\tilde {A}_K ]+[{W}_1 , \dots , {W}_K ]+ [\tilde {B}_1 ,
\dots ,\tilde {B}_K ] \}} H(\bar{t}){\rm d}\bar{t}
\]
where
\[
X_k = [x_{k0} , \dots ,x_{k,M - 1} ], \ X_{0k} = [x_0 , 0,\dots ,0
], \ k=1,2,\dots,K.
\]
Thus
\[
[X_1 , \dots ,X_K ] - [X_{01} , \dots ,X_{0K} ]= \{[  \tilde {A}_1
, \dots ,\tilde {A}_K ]+ [ {W}_1 , \dots ,{W}_K ]+ [\tilde {B}_1 ,
\dots ,\tilde {B}_K ] \}P.
\]
Using Kronecker product we rewrite the above equation as
\[
\hat {X} - \hat {X}_0 = (P^\tau \otimes I_n )\{[  \hat {\tilde
{A}}_1^\tau , \dots ,\hat {\tilde {A}}_K^\tau ]^\tau  +
[\hat{W}_1^\tau , \dots ,\hat{W}_K^\tau ]^\tau+[\hat {\tilde
{B}}_1^\tau , \dots ,\hat {\tilde {B}}_K^\tau ]^\tau \},
\]
where
\[
\hat {X} = [\hat {X}_1^\tau , \dots ,\hat {X}_K^\tau ]^\tau , \
\hat {X}_k = [x_{k0}^\tau , \dots ,x_{k,M - 1}^\tau ]^\tau , \
\hat {\tilde {A}}_k = [\tilde {A}_{k0}^\tau , \dots ,\tilde
{A}_{k,M - 1}^\tau ]^\tau,
\]
\[
 \hat {X}_0 = [\hat{X}_{01}^\tau , \dots
,\hat{X}_{0K}^\tau ]^\tau, \ \hat{X}_{0k} = [x^\tau_0 ,
0^\tau,\dots ,0^\tau ]^\tau, \ k=1,2,\dots,K.
\]
$\hat {W}_K$ and $\hat {\tilde {B}}_K$ have the similar meaning as
$\hat {\tilde {A}}_K$. So
\[
\hat {X} - \hat {X}_0 = (P^\tau \otimes I_n )\{[ (\hat {A}_1 \hat
{X}_1 )^\tau , \dots ,(\hat {A}_K \hat {X}_K )^\tau ]^\tau  +
\hat{W}
 + [(\hat {B}_1 \hat{U}_1 )^\tau , \dots
,(\hat {B}_K \hat{U}_K )^\tau ]^\tau\}
\]
\[
 \hspace*{40pt}= (P^\tau \otimes I_n )\bigg[\sum\limits_{k = 1}^K ({E_{kk}^{(K)} \otimes \hat
 {A}_k
} )\hat {X} +\sum\limits_{k = 1}^{[\frac{M+1}{2} ]} \sum\limits_{i
= 1}^K {\sum\limits_{j = 1}^{K}E^{(K)}_{ij}\otimes
\frac{d_j}{2k-1}{{N}^{(i)}_{j,2k-1}\hat{X}}}
\]
\[
\hspace*{192pt}  + \sum\limits_{k= 1}^K (E_{kk}^{(K)} \otimes \hat
{B}_k  )\hat {U}\bigg].
\]
Therefore
\begin{equation}
\hat {X}= \mathit\Gamma\hat {U}+ \mathit\Omega,\label{X}
\end{equation}
where
\[
\mathit\Gamma=[I_{MKn} - (P^\tau \otimes I_n
)\mathit\Phi]^{-1}(P^\tau\otimes I_n  )\sum\limits_{k= 1}^K
({E_{kk}^{(K)} \otimes \hat {B}_k } ),
\]
\[
\mathit\Phi=\sum\limits_{k = 1}^K (E_{kk}^{(K)} \otimes \hat
 {A}_k
 ) +\sum\limits_{k = 1}^{[\frac{M+1}{2} ]} \sum\limits_{i = 1}^K
\sum\limits_{j = 1}^{K}E^{(K)}_{ij}\otimes
\frac{d_j}{2k-1}{{N}^{(i)}_{j,2k-1}},
\]
\[
\mathit\Omega=[I_{MNn} - (P^\tau \otimes I_n
)\mathit\Phi]^{-1}\hat{X}_0.
\]
\\

\noindent\textbf{4. Numerical examples}
\\

\normalsize\noindent {\it 4.1. Example 1}
\\

Consider the integrodifferential system described by
\[
 \begin{array}{l}
\dot{x}(t)=\left[\begin{array}{cc}t^2+1&-t\\
0&1\end{array}\right]x(t)+\displaystyle\int_{\,0}^{\,1}\left[\begin{array}{cc}\bar{t}&3\\
3t^2&0\end{array}\right]x(\bar{t}){\rm d}\bar{t}+\left[\begin{array}{c}-(t-1)^2\\
2t^2-t^3\end{array}\right]u(t), \\
x(0)=[0,0]^\tau,
\end{array}
\]
where $x(t)$ is a 2-dimensional state function and $u(t)$ a
1-dimensional control function.

Let $u(t)=1$. Then by Eq.(\ref{X}) and choosing $M=4$ and $K=3$,
we have the hybrid solutions $x(t)=[x_1(t),x_2(t)]^\tau$:
\[
\begin{array}{l}
x_1(t)=\frac{5}{144}h_{10}(t)+\frac{1}{36}h_{11}(t)+\frac{1}{144}h_{12}(t)+\frac{37}{144}h_{20}(t)+\frac{1}{12}h_{21}(t)+\frac{1}{144}h_{22}(t)
\end{array}
\]
\[
\hspace*{30pt}
\begin{array}{l}
+\frac{101}{144}h_{30}(t)+\frac{5}{36}h_{31}(t)+\frac{1}{144}h_{32}(t),
\end{array}
\]
\[
\begin{array}{l}
x_2(t)=\frac{7}{864}h_{10}(t)+\frac{7}{864}h_{11}(t)+\frac{1}{128}h_{12}(t)
+\frac{1}{1728}h_{13}(t)+\frac{13}{96}h_{20}(t)+\frac{55}{864}h_{21}(t)
\end{array}
\]
\[
\hspace*{30pt}
\begin{array}{l}
+\frac{1}{96}h_{22}(t) +\frac{1}{1728}h_{23}(t)
+\frac{515}{864}h_{30}(t)+\frac{151}{864}h_{31}(t)+\frac{5}{288}h_{32}(t)
+\frac{1}{1728}h_{33}(t),
\end{array}
\]
which can be checked to be the same as the exact solutions
$x(t)=[t^2,t^3]^\tau$.
\\
\\

\normalsize\noindent {\it 4.2. Example 2}
\\

Consider the integrodifferential system described by
\[
 \begin{array}{l}
\dot{x}(t)=\left[\begin{array}{cc}1&t\\
t&t^2+1\end{array}\right]x(t)+\displaystyle\int_{\,0}^{\,1}\left[\begin{array}{cc}3\bar{t}&e^{-t}-\bar{t}^2\\
3t^2+\bar{t}e^{-t}&-t^2\end{array}\right]x(\bar{t}){\rm d}\bar{t}\\\hspace{130pt} +\left[\begin{array}{c}3e^{-1}-5-3t\\
2e^{-1}-7-t-3t^2\end{array}\right]u(t) , \ x(0)=[1,3]^\tau,
\end{array}
\]
where $x(t)$ is a 2-dimensional state function and $u(t)$ a
1-dimensional control function.

Let $u(t)=e^{-t}$. Then by Eq.(\ref{X}) and choosing $M=5$ and
$K=3$, we have the hybrid solutions $x(t)=[x_1(t),x_2(t)]^\tau$:
\[
\begin{array}{l}
x_1(t)=\frac{959}{1129}h_{10}(t)-\frac{202}{2857}h_{11}(t)+\frac{67}{22756}h_{12}(t)-\frac{3}{36694}h_{13}(t)+\frac{1}{587104}h_{14}(t)
\end{array}
\]
\[
\hspace*{30pt}
\begin{array}{l}
+\frac{465}{764}h_{20}(t)-\frac{203}{4007}h_{21}(t)+\frac{121}{57355}h_{22}(t)-\frac{5}{85351}h_{23}(t)+\frac{1}{819370}h_{24}(t)
\end{array}
\]
\[
\hspace*{30pt}
\begin{array}{l}
+\frac{529}{1213}h_{30}(t)-\frac{261}{7190}h_{31}(t)+\frac{47}{31092}h_{32}(t)-\frac{5}{119117}h_{33}(t),
\end{array}
\]
\[
\begin{array}{l}
x_2(t)=\frac{2877}{1129}h_{10}(t)-\frac{613}{2890}h_{11}(t)+\frac{187}{21171}h_{12}(t)
-\frac{9}{36694}h_{13}(t)+\frac{2}{391403}h_{14}(t)
\end{array}
\]
\[
\hspace*{30pt}
\begin{array}{l}
+\frac{944}{517}h_{20}(t)-\frac{226}{1487}h_{21}(t)+\frac{364}{57513}h_{22}(t)
-\frac{14}{79661}h_{23}(t)+\frac{2}{546247}h_{24}(t)
\end{array}
\]
\[
\hspace*{30pt}
\begin{array}{l}
+\frac{1744}{1333}h_{30}(t)-\frac{449}{4123}h_{31}(t)+\frac{47}{10364}h_{32}(t)
-\frac{8}{63529}h_{33}(t)+\frac{1}{381175}h_{34}(t).
\end{array}
\]
\begin{table}[htbp]
  \caption{\small Analytic and hybrid values of  $x_1(t)$.}
\begin{tabular}{|c|c|c|c|c|} \hline
 \small  $t$ &\small   Analytic & \small   Hybrid ($N$=4,$M$=5)&\small Hybrid ($N$=4,$M$=7)& \small Hybrid ($N$=4,$M$=9)\\\hline
\small0.1&\small0.90483741803596 &\small0.90483741135846 &\small0.90483741803662&\small0.90483741803596\\
\small0.2&\small0.81873075307798 &\small0.81873074696139 &\small0.81873075307870&\small0.81873075307798\\
\small0.3&\small0.74081822068172 &\small0.74081822463206 &\small0.74081822068111&\small0.74081822068172\\
\small0.4&\small0.67032004603564 &\small0.67032005059525 &\small0.67032004603507&\small0.67032004603564\\
\small0.5&\small0.60653065971263 &\small0.60653063311834 &\small0.60653065970940&\small0.60653065971263\\
\small0.6&\small0.54881163609403 &\small0.54881163079527 &\small0.54881163609435&\small0.54881163609402\\
\small0.7&\small0.49658530379141 &\small0.49658529860747 &\small0.49658530379175&\small0.49658530379140\\
\small0.8&\small0.44932896411722 &\small0.44932896475352 &\small0.44932896411673&\small0.44932896411722\\
\small0.9&\small0.40656965974060 &\small0.40656966036876 &\small0.40656965974011&\small0.40656965974059\\
\small1.0&\small0.36787944117144 &\small0.36787945775656&\small 0.36787944117369&\small0.36787944117143\\
\hline
\end{tabular}
\end{table}
\begin{table}[htbp]
  \caption{\small Analytic and hybrid values of  $x_2(t)$.}
\begin{tabular}{|c|c|c|c|c|} \hline
 \small  $t$ &\small   Analytic & \small   Hybrid ($N$=4,$M$=5)&\small Hybrid ($N$=4,$M$=7)& \small Hybrid ($N$=4,$M$=9)\\\hline
\small0.1&\small2.71451225410788 &\small2.71451223459708 &\small2.71451225410991&\small2.71451225410788\\
\small0.2&\small2.45619225923395 &\small2.45619224196066 &\small2.45619225923618&\small2.45619225923394\\
\small0.3&\small2.22245466204515 &\small2.22245467555371 &\small2.22245466204343&\small2.22245466204515\\
\small0.4&\small2.01096013810692 &\small2.01096015405462 &\small2.01096013810538&\small2.01096013810692\\
\small0.5&\small1.81959197913790 &\small1.81959190225067 &\small1.81959197912838&\small1.81959197913790\\
\small0.6&\small1.64643490828208 &\small1.64643489597035 &\small1.64643490828328&\small1.64643490828208\\
\small0.7&\small1.48975591137423 &\small1.48975590012432 &\small1.48975591137553&\small1.48975591137423\\
\small0.8&\small1.34798689235166 &\small1.34798689933842 &\small1.34798689235054&\small1.34798689235166\\
\small0.9&\small1.21970897922180 &\small1.21970898705545 &\small1.21970897922074&\small1.21970897922179\\
\small1.0&\small1.10363832351433 &\small1.10363838023956 &\small1.10363832352153&\small1.10363832351432\\
\hline
\end{tabular}
\end{table}
Table 1 and Table 2 imply that  the hybrid solutions converges the
exact solutions $x(t)=[e^{-t},3e^{-t}]^\tau$ rapidly.
\\

\noindent\textbf{6. Conclusion}

Using the excellent properties of operational matrices of the
hybrid function of general block-pulse functions and the second
Chebyshev polynomials, the general algorithms for
integrodifferential systems are derived. The illustrative examples
demonstrate that this technique is convenient for application.
\\

\noindent{\bf Acknowledgement}
\\

This work was supported by Science Research Foundation in Harbin
Institute of Technology (Grant No. HITC200708).

\def\refname

{\normalsize\textbf References}
\begin{thebibliography}{1}
\footnotesize

\bibitem{Kochetkov}Yu.A. Kochetkov, V.P. Tomshin, Optimal control of deterministic system described by integrodifferential equations.
 Automat. Remote Control 39 (1978) 1-6.

\bibitem{block1}P. Sannuti, Analysis and synthesis of dynamic systems via block-pulse functions.
 Proc. Inst. Elect. Eng. 124 (1977) 569-571.

\bibitem{block2}G.P. Rao, L. Sivakumar, Analysis and synthesis of dynamic systems containing
time-delays via block-pulse functions.
 Proc. IEE 125 (1978) 1064-1068.


\bibitem{Walsh}K.R. Palanisamy, G.P. Rao, Optimal control of
linear systems with delays in state and control via Walsh
function.  Proc. Inst. Elect. Eng. 130 (1983) 300-312.

\bibitem{Fourier}M. Razzaghi, M. Razzaghi, Fourier series approach for the solution of Linear
two-point boundary value problems with time-varying coefficients.
 Int. J. Syst. Sci.
21 (1990) 1783-1794.

\bibitem{Legendre}H. Lee, F.C. Kung, Shifted Legendre series solution
and parameter estimation of linear delayed systems.
 Int. J.
Syst. Sci. 16 (1985) 1249-1256.

\bibitem{Chebyshev}I.R. Horng, J.H. Chou, Analysis, parameter
estimation and optimal control of time-delay systems via Chebyshev
series.  Int. J. Control 41 (1985) 1221-1234.

\bibitem{Laguerre}C. Hwang, Y.P. Shih,  Laguerre series solution of a functional differential equation,   Int. J. Syst. Sci. 13 (1982) 783--788.

\bibitem{hybridTaylor}K. Maleknejad, Y. Mahmoudi, Numerical solution of
linear Fredholm integral equation by using hybrid Taylor and
block-pulse functions. Appl. Math. Comput. 149 (2004) 799-806.


\bibitem{hybridLegendre}X.T. Wang, Numerical solutions of optimal control for time delay systems  by hybrid of block-pulse functions and
Legendre polynomials,  Appl. Math. Comput. 184 (2007) 849--856.

\bibitem{Legendrewavelets}X.T. Wang, Numerical solution of  time-varying systems with a stretch by
general Legendre wavelets,  Appl. Math. Comput. 198 (2008)
613--620.

\end{thebibliography}
\end{document}